\documentclass[12pt]{amsart}
\usepackage{latexsym}
\usepackage{amsfonts}
\usepackage{amsthm}
\usepackage{amsmath}
\usepackage{amssymb}

\def\sw#1{{\sb{(#1)}}}

\def\<{{\langle}}
\def\>{{\rangle}}

\def\eps{\epsilon}

\def\note#1{{}}
\def\can{{\rm can}}

\def\note#1{}

\def\Z{{\mathbb{Z}}}

\def\can{{\rm can}}

\def\beq{\begin{equation}}
\def\eeq{\end{equation}}

\def\id{{I}}

\def\ut{{\otimes}}
\def\ot{{\otimes}}

\newcounter{zlist}

\def\Label#1{\label{#1}\ifmmode\llap{[#1] }\else
\marginpar{\smash{\hbox{\tiny [#1]}}}\fi}
\def\Label{\label}

\newtheorem{proposition}{Proposition}[section]

\newtheorem{theorem}[proposition]{Theorem}

\theoremstyle{definition}
\newtheorem{definition}[proposition]{Definition}
\newtheorem{example}[proposition]{Example}

\theoremstyle{remark}
\newtheorem{remark}[proposition]{Remark}

\newcounter{c}

\newcommand{\etyk}[1]{\vspace{-7.4mm}$$\begin{equation}\Label{#1}
\addtocounter{c}{1}}
\renewcommand{\]}{\ifnum \value{c}=1 $$\else \end{equation}\fi}
\begin{document}

\title{Yang-Baxter systems and
Entwining structures}
\author{Tomasz Brzezi\'nski}
\address{Department of Mathematics, University of Wales Swansea,
Singleton Park, Swansea SA2 8PP, U.K.} 
\email{T.Brzezinski@swansea.ac.uk}
\author{Florin F. Nichita}
\address{Institute of Mathematics of the Romanian Academy, 
P.O. Box 1-764, RO-70700, Bucharest, Romania}
\address{Department of Mathematics, University of Wales Swansea,
Singleton Park, Swansea SA2 8PP, U.K.}
\email{F.F.Nichita@swansea.ac.uk}
\date{November 2003}
\subjclass{16W30, 17B37, 81R50}
\begin{abstract}
It is shown that a Yang-Baxter system can be constructed from 
any entwining structure. It is also shown that, conversely, Yang-Baxter 
systems of certain type lead to entwining structures. Examples of Yang-Baxter 
systems associated to entwining structures are given, and 
a Yang-Baxter operator of Hecke type is defined for any bijective entwining map.
\end{abstract}
\maketitle
\section{Introduction and preliminaries}

The origins of the quantum Yang-Baxter equation are  in theoretical physics, where it first appeared in the context of factorisation of $S$-matrices and  in the quantum inverse scattering method of solving integrable Hamiltonian systems. Since the work of Drinfeld and Jimbo, Yang-Baxter equations are recognised as  the origin of the theory of quantum groups or quasitriangular Hopf algebras. They are the starting point of many developments in knot theory and  in  the theory of braided monoidal categories. Due to the complexity of various integrable models the need arose for extensions and generalisations of Yang-Baxter equations and related algebraic structures.  One such  generalisation, termed a {\em Yang-Baxter system}, was proposed in  \cite{HlaSno:sol} in the context of non-ultralocal integrable systems  \cite{HlaKun:qua}. 

Motivated by the need for developing a general theory of non-commutative principal bundles, entwining structures were introduced in \cite{BrzMaj:coa} as generalised symmetries of such bundles. In recent years, entwining structures and categories of modules associated to them have been used to provide a unification of various types of Hopf modules (cf.\ \cite{Brz:mod}, \cite{CaeMil:Fro}), and eventually have led to the revival of the theory of corings (cf.\ \cite{Brz:str}, \cite{BrzWis:cor}).

Thus Yang-Baxter systems and entwining structures arose in entirely different contexts, were motivated by as far fields as high energy physics and non-commutative geometry, and have entirely different applications. There is no reason to expect that there is any connection between these two notions. Yet, in this paper, we prove that such a connection, and a very close one at that, indeed exists. More precisely, we show that to any entwining
structure one can associate a Yang-Baxter system. The converse  is also
true: to any Yang-Baxter
system of certain type one can associate an
entwining structure.

The remainder of this section contains preliminaries about Yang-Baxter operators. The main results and examples are contained in Section~2. We first recall the notions of a Yang-Baxter system and an entwining structure and then prove the theorem connecting these two notions. We then present a series of examples of Yang-Baxter systems related to various entwining structures. We also use   the classification of two-dimensional Yang-Baxter systems
given in \cite{HlaSno:sol} to construct an explicit 
example of an entwining structure. Finally, we show that to every 
entwining structure $(A,C)_\psi$ with a bijective entwining map $\psi$, 
one can associate a $q$-Hecke operator on $A\oplus C$ obtained by gluing $q$-Hecke operators on $A$ and $C$. 

Throughout this paper $k$ is a field. The unadorned tensor product is over $k$. The identity map on a $k$-vector space $V$ is denoted by $\id_V$ or simply by $\id$ provided the domain is clear from the context. All algebras are over $k$, 
they are associative and with unit 1.  
The product in an algebra $A$ is denoted by $\mu:A\ot A\to A$, while the unit map is denoted by $\iota:k\to A$. Similarly, all coalgebras are over $k$, they are coassociative  and with a counit.  Coproduct in a coalgebra $C$ is denoted by $\Delta:C\to C\ot C$ and the counit by  $\eps:C\to k$.   
We use the standard notation for coalgebras and comodules. In particular, for a coalgebra $C$, we use Sweedler's notation to denote the coproduct $\Delta$ on elements, i.e., $\Delta(c) = \sum c\sw1\ot c\sw 2$, for all $c\in C$. For a right $C$-comodule $V$, the coaction is denoted by $\varrho^V(v) = \sum  v\sw 0\ot v\sw 1$, for all $v\in V$.

For any vector spaces $ V $ and $ W $,
$ \tau_{V, W}: V \otimes W \rightarrow W \ot V \  $ denotes the natural
bijection defined by $ \tau_{V, W}(v \ot w) = w \ot v $.
Let $ R: V \ot V \rightarrow V \ot V  $
be a $ k$-linear map. Define
$ {R_{12}}= R \ot \id_V$, ${R_{23}}= \id_V \ot R$  and
${R_{13}}=(\id_V\ot\tau_{V, V})\circ(R\ot \id_V)\circ (\id_V\ot \tau_{V, V})$. Each of the $R_{ij}$ is thus a linear endomorphism of $V\ot V\ot V$. Recall that an invertible  $k$-linear map  $ R : V \ot V \rightarrow V \ot V $
is called a {\em Yang-Baxter
operator} (or simply a {\em YB operator}) if it satisfies the  equation
\begin{equation}  \label{ybeq}
R_{12}  \circ  R_{23}  \circ  R_{12} = R_{23}  \circ  R_{12}  \circ  R_{23}.
\end{equation}
The equation (\ref{ybeq}) is usually called the {\em braid equation}. It is a
well-known fact that the operator $R$ satisfies (\ref{ybeq}) if and only if
$R\circ \tau_{V, V}  $ satisfies
   the {\em quantum Yang-Baxter equation} 
\begin{equation}   \label{ybeq2}
R_{12}  \circ  R_{13}  \circ  R_{23} = R_{23}  \circ  R_{13}  \circ  R_{12},
\end{equation}
perhaps more familiar in physics and quantum group theory. For a review of Yang-Baxter operators we refer to \cite{lara}.

In what follows, we use the following construction of Yang-Baxter operators described in \cite{DasNic:yan}. If $A$ is a $k$-algebra, then for all non-zero $r,s\in k$, the linear map 
\begin{equation}\label{ra}
R^A_{r,s}:A \ot A
\rightarrow A \ot A,\quad   a \ot b \mapsto s ab \ot 1 + r 1 \ot ab - s a \ot b 
\end{equation}
is a Yang-Baxter operator. The inverse of $R^A_{r,s}$ is 
${(R^A_{r,s})}^{-1}(a \ot b)= \frac{1}{r} ab \ot 1 + \frac{1}{s} 1 \ot 
ab - \frac{1}{s} a \ot b $.
Dually, if $C$ is a coalgebra, then for all non-zero $p,t\in k$, the linear map 
\begin{equation}\label{rc}
R_C^{p,t}:C \ot C
\rightarrow C \ot C,\quad   c \ot d \mapsto p \eps (c) \sum d\sw 1 \ot d\sw 2 + t \eps (d) 
\sum c\sw 1 \ot c\sw 2 - p c \ot d 
\end{equation}
is a Yang-Baxter operator. Note that in both cases the assumption that parameters are non-zero is needed only for the invertibility of $R^A_{r,s}$ and $R_C^{p,t}$. $R^A_{r,s}$ and $R_C^{p,t}$ satisfy the braid relation for any value of $r,s,p$ and $t$.

\section{From entwining structures to Yang-Baxter systems (and back)}

Yang-Baxter systems were introduced in \cite{HlaSno:sol} as a spectral-parameter independent generalisation of quantum Yang-Baxter equations related to non-ultralocal integrable systems studied previously in  \cite{HlaKun:qua}. Yang-Baxter systems are conveniently defined in terms of {\em Yang-Baxter commutators}. Consider three vector spaces $V,V',V''$ and  three linear maps
$ R : V \ot V' \rightarrow V \ot V' $,
$ S : V \ot V'' \rightarrow V \ot V'' $ and
$ T : V' \ot V'' \rightarrow V' \ot V'' $. Then a {\em Yang-Baxter commutator} is a map
$ [R,S,T]:  V \ot V' \ot V'' \rightarrow V \ot V' \ot V'' $, defined by
\begin{equation}   \label{ybcomm}
[R,S,T]= R_{12}  \circ  S_{13}  \circ  T_{23} - T_{23}  \circ  S_{13} 
\circ  R_{12} \ .
\end{equation}
In terms of a Yang-Baxter commutator,  the quantum Yang-Baxter equation (\ref{ybeq2})  is expressed simply as $[R,R,R] = 0$.

\begin{definition}\label{def.wxz}
Let $V$ and  $V'$ be vector spaces. A system of linear maps 
$$ W : V \ot V \rightarrow V \ot V ,\quad 
Z : V' \ot V' \rightarrow V' \ot V' , \quad
 X : V \ot V' \rightarrow V \ot V' 
$$ 
is called  a {\em WXZ-system} or a 
{\em Yang-Baxter system}, provided the following equations are satisified:
\begin{equation}   \label{ybeqn4}
[W,W,W]\ = \  0 \ ,
\end{equation}
\begin{equation}   \label{ybeqn5}
[Z,Z,Z]\ = \  0 \ ,
\end{equation}
\begin{equation}   \label{ybeqn6}
[W,X,X]\ = \  0 \ ,
\end{equation}
\begin{equation}   \label{ybeqn7}
[X,X,Z]\ = \  0 \ .
\end{equation}
\end{definition}
\bigskip

There are several algebraic origins and applications of WXZ-systems. It has been observed in \cite{Vla:met} that WXZ-systems with invertible $W$, $X$ and $Z$ can be used to construct dually-paired bialgebras of the FRT type, thus leading to quantum doubles. More precisely, consider a WXZ-system with finite-dimensional $V=V'$, so that each of $W$, $X$, $Z$ is an $N^2\times N^2$-matrix. Suppose that $W,X,Z$ are invertible. Since $W$ and $Z$ satisfy Yang-Baxter equations (\ref{ybeqn4})--(\ref{ybeqn5}), one can consider two matrix bialgebras $A$ and $B$ with $N\times N$ matrices of generators $U$ and $T$ respectively, and relations
$ W_{12} U_{1} U_{2} = U_{2} U_{1} W_{12} $, $ Z_{12} T_{1} T_{2} = T_{2} T_{1} Z_{12}$ (cf.\ \cite{FadRes:Lie}). 
The existence of 
an invertible operator $X$ that satisfies equations (\ref{ybeqn6})--(\ref{ybeqn7}), means that $A$ and $B$ are dually paired with a non-degenerate pairing $ <U_{1},T_{2}> = X_{12} $. Furthermore,  the tensor product $A\ot B$ has an algebra (quantum double) structure with crossed relations
$ X_{12}U_{1}T_{2} = T_{2}U_{1}X_{12} $.

Given a WXZ-system as in Definition~\ref{def.wxz} one can construct a Yang-Baxter operator on $V\oplus V'$, provided the map $X$ is invertible. This is a special case of a {\em gluing procedure} described in  \cite[Theorem~2.7]{MajMar:glu} (cf.\  \cite[Example 2.11]{MajMar:glu}).
Let $ R = W \circ \tau_{V,V} $, $ R' = Z \circ \tau_{V',V'} $, $ U = X
\circ \tau_{V',V} $. Then the linear map
$$
 R \oplus_{U} R': (V \oplus V') \ot (V \oplus V') \to  (V \oplus V') \ot (V \oplus V')
$$ 
given by
$ R \oplus_{U} R'|_{V \ot V} = R $,
$ R \oplus_{U} R'|_{V' \ot V'} = R' $, and for all $ x \in V $, $ y \in V' $,
$$ 
(R \oplus_{U} R')(y \ot x) = U(y \ot x ), \qquad  (R \oplus_{U} R')(x \ot y ) = U^{-1}(x \ot y)
$$
is a Yang-Baxter operator.

Entwining structures were introduced in \cite{BrzMaj:coa} in order to recapture the symmetry structure of non-commutative (coalgebra) principal bundles or coalgebra-Galois extensions. For  applications and algebraic content of entwining structures we refer to  \cite{CaeMil:Fro} and \cite{BrzWis:cor}.

\begin{definition}\label{def.entw} 
 An algebra $A$ is said to be {\em entwined} with  a coalgebra $C$ if there exists  a linear  map 
$ \psi : C \ot A
\rightarrow A \ot C $ satisfying the following four conditions:

(1) $ \psi \circ (I_{C} \ot \mu) = ( \mu \ot I_{C}) \circ (I_{A} \ot
\psi) \circ ( \psi \ot I_{A}) $,

(2) $ (I_{A} \ot \Delta) \circ \psi = ( \psi \ot I_{C} ) \circ (I_{C}
\ot \psi) \circ ( \Delta \ot I_{A} ) $,

(3) $ \psi \circ ( I_{C} \ot \iota ) = \iota \ot I_{C} $,

(4) $ (I_{A} \ot \eps ) \circ \psi = \eps \ot I_{A} $.

The map $ \psi $ is known as an {\em entwining map}, and the triple  $ {(A,C)}_{\psi} $ is called an {\em entwining structure}. 
\end{definition}

To denote the action of an entwining
map $ \psi $ on elements it is convenient to use the following {\em $\alpha$-notation}, for all $a,b\in A$ and $c\in C$,
$$
\psi (c \ot a) = \sum_{\alpha} a_{\alpha} \ot c^{\alpha} , \quad
(I_{A} \ot \psi) \circ ( \psi \ot I_{A}) (c \ot a \ot b) =
\sum_{\alpha, \beta} a_{\alpha} \ot b_{\beta} \ot c^{\alpha \beta},
$$
etc. The relations (1), (2), (3)
and (4) in Definition~\ref{def.entw} are  equivalent to the following explicit relations, for all $ a,b
\in A$, $c \in C$,
\begin{equation}\label{a}
\sum_{\alpha}(ab)_{\alpha}\ot c^{\alpha} = \sum_{\alpha, \beta}
a_{\alpha} b_{\beta} \ot c^{\alpha \beta}
\end{equation}
\begin{equation}\label{b}
\sum_{\alpha} a_{\alpha} \ot {c^{\alpha}}\sw{1} \ot {c^{\alpha}}\sw{2} =
\sum_{\alpha, \beta} a_{\beta \alpha} \ot {c\sw{1}}^{\alpha} \ot
{c\sw{2}}^{\beta}
\end{equation}
\begin{equation}\label{c}
\sum_{\alpha} 1_{\alpha} \ot c^{\alpha} = 1 \ot c
\end{equation}
\begin{equation}\label{d}
\sum_{\alpha}a_{\alpha} \eps (c^{\alpha}) = a \eps (c)
\end{equation}

The main result of this paper is contained in the following
\begin{theorem}\label{thm.main}
 Let $A$ be an algebra  and let $C$ be a
coalgebra. For any $ s, r, t, p \in k$ define linear maps 
$$W: A\ot A\to A\ot A, \qquad a\ot b \mapsto s ba \ot 1 + r 1 \ot ba - s b \ot a, 
$$
$$
Z: C\ot C\to C\ot C, \qquad c\ot d \mapsto t \eps (c) \sum d\sw 1 \ot d\sw 2 + p \eps (d) 
\sum c\sw 1 \ot c\sw 2 - p d \ot c.
$$
Let $X:A\ot C\to A\ot C$ be a linear map such that $ X\circ (\iota \ot \id_C) = \iota\ot \id_{C} $ and $ (\id_A \ot \eps ) \circ X =\id_A\ot \eps$. Then $W,X,Z$ is a Yang-Baxter system if and only if $A$ is entwined with $C$ by the map $\psi: = X\circ \tau_{C, A}$.

\end{theorem}

\begin{proof}
 Note that $W = R^A_{r,s}\circ\tau_{A,A}$ and $Z= R_C^{p,t}\circ\tau_{C,C}$, where $R^A_{r,s}$ is defined  in (\ref{ra}), and  $ R_C^{p,t}$ is defined in (\ref{rc}), so that the conditions $[W,W,W]=0$ and $[Z,Z,Z]=0$ are satisfied. 

Suppose that $(A,C)_\psi$ is an entwining structure. Equations (\ref{c})--(\ref{d}) imply that $ X\circ (\iota \ot \id_C) = \iota\ot \id_{C} $ and $ (\id_A \ot \eps ) \circ X =\id_A\ot \eps $. The fact that $X = \psi\circ\tau_{A,C}$ satisfies conditions (\ref{ybeqn6})--(\ref{ybeqn7}) can be checked by direct computations. Explicitly, using the $\alpha$-notation we obtain $X(a\ot c) = \sum_{\alpha} a_{\alpha}\ot c^{\alpha}$, for all $a\in A$ and $c\in C$. In view of the definition of $W$, 
$W_{12} (a \ot b \ot c) = s ba \ot 1 \ot c + r 1 \ot ba \ot c -
s b \ot a \ot c $, for all $a,b\in A$ and $c\in C$, so that
$$
X_{13} \circ W_{12} (a \ot b \ot c ) =
s \sum_{\alpha} (ba)_{\alpha} \ot 1 \ot c^{\alpha} + r 1 \ot ba \ot c -
s \sum_{\alpha} b_{\alpha} \ot a \ot c^{\alpha} .
$$
The second expression on the right-hand side is obtained with the help of formula (\ref{c}). Again the use of  (\ref{c}) leads to
$$ 
X_{23} \circ X_{13} \circ W_{12} (a \ot b \ot c)= s \sum_{\alpha}
(ba)_{\alpha} \ot 1 \ot c^{\alpha} + r \sum_{\alpha} 1 \ot 
(ba)_{\alpha} \ot c^{\alpha} -
s \sum_{\alpha , \beta } b_{\alpha} \ot a_{\beta} \ot c^{\alpha \beta}.
$$
On the other hand, the definition of $X$ yields
$X_{23} (a \ot b \ot c) = \sum_{\alpha} a \ot b_{\alpha} \ot
c^{\alpha}$, so that
$X_{13} \circ X_{23} (a \ot b \ot c)=\sum_{\alpha , \beta} a_{\beta} 
\ot b_{\alpha} \ot
c^{\alpha \beta}$, and therefore
$$ 
W_{12} \circ X_{13} \circ X_{23} (a \ot b \ot c)= s \sum_{\alpha, \beta}
b_{\alpha} a_{\beta}
  \ot 1 \ot c^{\alpha \beta} + r \sum_{\alpha,\beta} 1 \ot  b_{\alpha} a_{\beta}
\ot c^{\alpha \beta} - s \sum_{\alpha , \beta } b_{\alpha} \ot 
a_{\beta} \ot c^{\alpha \beta}.
$$
Now, using equation  (\ref{a}) for the first and second terms on the right-hand side, we obtain
$ [ W, X, X] = 0$, as required. Thus  (\ref{ybeqn6}) holds.

The equality $ [ X, X, Z ]=0 $ can be deduced from the self-duality of an entwining structure. More precisely, the conditions Definition~\ref{def.entw} are invariant under the operation consisting of exchanging $A$ with $C$, $\mu$ with $\Delta$,  $\iota$ with $\eps$ and reversing the order of the composition (no change in $\psi$). As a result of this operation, condition (1) in Definition~\ref{def.entw} becomes condition (2), condition (2) becomes condition (1), and (3) is exchanged with (4) - but there is no overall change in the content of Definition~\ref{def.entw}.  The above operation changes $W$ to a map of the form $Z$, and the Yang-Baxter commutator $[ W, X, X] $ to the commutator of the form  $[ X, X, Z]$ . Since the calculation of $[ W, X, X] $ does not depend on the choice of the parameters $r$, $s$, $[ W, X, X] =0$ implies $[ X, X, Z ]=0$. This can also be checked explicitly. Directly from the definitions of $X$ and $Z$ one finds, for all $a\in A$ and $c,d\in C$,

\begin{eqnarray*}
X_{12} \circ X_{13} \circ Z_{23} (a \ot c \ot d) &=&
t \eps (c) \sum_{\beta , \alpha} a_{\beta \alpha} \ot {d\sw 1}^{\alpha} \ot
{d\sw 2}^{\beta } \\
&&+ p \eps (d) \sum_{\beta , \alpha} a_{\beta \alpha} 
\ot {c\sw 1}^{\alpha} \ot {c\sw 2}^{\beta } - p \sum_{ \beta , \alpha} a_{\beta \alpha} \ot
{d}^{\alpha} \ot {c}^{\beta} 
\end{eqnarray*}
and

\begin{eqnarray*}
Z_{23} \circ X_{13} \circ X_{12} (a \ot c \ot d) &=&
t \sum_{\alpha, \beta}  \eps (c^{\alpha})
a_{\alpha \beta} \ot {d}^{\beta}\sw 1 \ot {d}^{\beta }\sw 2\\
&&+
p \sum_{\alpha, \beta} a_{\alpha \beta}
\eps (d^{\beta})
  \ot {{c}^{\alpha}}\sw 1 \ot {{c}^{\alpha }}\sw 2
  - p \sum_{ \alpha, \beta} a_{\alpha \beta} \ot
{d}^{\beta} \ot {c}^{\alpha}.
\end{eqnarray*}
Next, apply equation (\ref{d}) (dual to equation (\ref{c})) to the first two terms on the right-hand side of the last equality to obtain

\begin{eqnarray*}
Z_{23} \circ X_{13} \circ X_{12} (a \ot c \ot d) &=&
t \sum_{\alpha, \beta}  \eps (c)
a_{\beta} \ot {d}^{\beta}\sw 1 \ot {d}^{\beta }\sw 2\\
&&+
p \sum_{\alpha} a_{\alpha}
\eps (d)
  \ot {{c}^{\alpha}}\sw 1 \ot {{c}^{\alpha }}\sw 2
  - p \sum_{ \alpha, \beta} a_{\alpha \beta} \ot
{d}^{\beta} \ot {c}^{\alpha}.
\end{eqnarray*}
Finally, application of equation (\ref{b}) (dual to equation (\ref{a})) to the first two terms on the right-hand side yields
 $[ X, X, Z ]=0$, as expected.

In the converse direction, suppose that the map $X:A\ot C\to A\ot C$ is such that $W$, $X$ and $Z$ form  a Yang-Baxter system
and that $ X\circ (\iota \ot \id_C) = \iota\ot \id_{C} $ and $ (\id_A \ot \eps ) \circ X =\id_A\ot \eps$. The latter immediately imply that $\psi = X\circ\tau_{C,A}$ satisfies conditions (3) and (4) in Definition~\ref{def.entw}. Write $\psi(c\ot a) = \sum_\alpha a_\alpha\ot c^\alpha$, and perform the same computations of $[W,X,X]$ and $[X,X,Z]$ as before (note that one can use equations (\ref{c})--(\ref{d})), to conclude that $[W,X,X]=0$ imply that, for all $a,b\in A$, $c\in C$,
$$
 \sum_{\alpha} 1 \ot 
(ba)_{\alpha} \ot c^{\alpha} = \sum_{\alpha,\beta} 1 \ot  b_{\alpha} a_{\beta}
\ot c^{\alpha \beta},
$$
so that equation (\ref{a}) holds. Furthermore, the equality $[X,X,Z]=0$ implies that, for all $a\in A$, $c,d\in C$,
$$
\sum_{\alpha} a_{\alpha}
\eps (d)
  \ot {{c}^{\alpha}}\sw 1 \ot {{c}^{\alpha }}\sw 2 = \eps (d) \sum_{\beta , \alpha} a_{\beta 
\alpha} \ot {c\sw 1}^{\alpha} \ot {c\sw 2}^{\beta } .
$$
It therefore follows that for all $a\in A$, $c,d\in C$,
$$
d\ot \sum_{\alpha} a_{\alpha}
  \ot {{c}^{\alpha}}\sw 1 \ot {{c}^{\alpha }}\sw 2 = d \ot \sum_{\beta , \alpha} a_{\beta 
\alpha} \ot {c\sw 1}^{\alpha} \ot {c\sw 2}^{\beta } .
$$
Since $k$ is a field, the above equality yields equation (\ref{b}). Thus $\psi$ is an entwining map.
\end{proof}

\begin{remark}
As observed in \cite[Proposition~2.7]{BrzMaj:coa}, entwinning structures are related to {\em algebra factorisations} via semi-dualisation (see \cite[p.\ 300]{Maj:fou} for discussion of algebra factorisations). The arguments similar to those in the proof of Theorem~\ref{thm.main} show that, given two algebras $A$, $B$, and an algebra factorisation map $\Psi: B\ot A\to A\ot B$, one can construct a Yang-Baxter system with $W$ and $Z$ of the same form as $W$ in Theorem~\ref{thm.main} and $X=\Psi\circ\tau_{A,B}$. 
\end{remark}

Thus Theorem~\ref{thm.main} implies that for any entwining structure $(A,C)_\psi$ there is a family (parametrised by fours elements of $k$) of Yang-Baxter systems. We describe two general
 examples below.
\begin{example} {\bf Yang-Baxter systems associated to a Doi-Koppinen datum}.
Let $B$ be a bialgebra, $A$ be a right B-comodule algebra and let  $C$ be a right $B$-module
coalgebra. Then, $A$ is entwined with $C$ via $ \psi : C \ot A \rightarrow A \ot C $, $ \psi (c \ot
a) = \sum a_{(0)} \ot c a_{(1)}$ (cf.\ \cite[Example~3.1]{Brz:mod}). Therefore, for all $r,s,p,t\in k$, the maps
$$
W: A\ot A\to A\ot A, \qquad a \ot b \mapsto  s ba \ot 1 + r 1 \ot ba - s b \ot a ,
$$
$$
X: A\ot C\to A\ot C, \qquad a \ot c\mapsto \sum  a_{(0)}\ot c a_{(1)} ,
$$
$$ Z:C\ot C\to C\ot C, \qquad c \ot d\mapsto  t \eps (c) \sum d\sw 1 \ot d\sw 2 +
p \eps (d) \sum c\sw 1 \ot c\sw 2 - p d \ot c 
$$
form a Yang-Baxter system. Furthermore, if $ B^{op} $ is a Hopf algebra
(i.e., there exists a map $ S: B \rightarrow B $ such that
$ \sum S(b_{(2)}) b_{(1)} = \sum b_{(2)} S(b_{(1)}) = \eps (b) $) then $ \psi $
is invertible and thus so is $X$.

As a particular case we can take $ A=C=B$
and view $B$ as a $B$-comodule algebra via $ \Delta $ and as a $B$-module
coalgebra via $\mu$. Thus any bialgebra gives rise to a Yang-Baxter system.
\end{example}

\begin{example} {\bf Yang-Baxter systems associated to coalgebra-Galois extensions}.
Let $C$ be a coalgebra, and $A$ be an algebra and a right $C$-comodule with a coaction $\rho^A:A\to A\ot C$. Define a subalgebra $B$ of $A$ by 
$$B = \{b\in A\; |\; \forall a\in A, \; \rho^A(ba) = b\rho^A(a) = \sum ba\sw 0\ot a\sw 1\}.
$$
 Suppose that the canonical right $C$-comodule, left $B$-module map 
$$\can: A\ot_B A\to A\ot C, \qquad a\ot a'\mapsto a\rho^A(a')$$ 
is bijective, i.e., that $B\subseteq A$ is a {\em coalgebra-Galois $C$-extension}. 
Then $A$ is entwined with $C$ via the map 
$\psi: C\ot A\to A\ot C$, $c\ot a\mapsto \can(\can^{-1}(1\ot c)a)$ 
(cf.\ \cite[Theorem~2.7]{BrzHaj:coa}). Hence, for all $r,s,p,t\in k$, the maps
$$
W: A\ot A\to A\ot A, \qquad a \ot b \mapsto  s ba \ot 1 + r 1 \ot ba - s b \ot a ,
$$
$$
X: A\ot C\to A\ot C, \qquad a \ot c\mapsto\can(\can^{-1}(1\ot c)a),
$$
$$ Z:C\ot C\to C\ot C, \qquad c \ot d\mapsto  t \eps (c) \sum d\sw 1 \ot d\sw 2 +
p \eps (d) \sum c\sw 1 \ot c\sw 2 - p d \ot c 
$$
form a Yang-Baxter system. 
\end{example}

The following example is based on the construction of an entwining structure in \cite[Example~3.4]{Sch:Doi}
\begin{example}
Let
 $C$ be a coalgebra spanned as a vector space  by $e$ and
      $y_j$ , $j\in \Z$
 with the coproduct $\Delta(e) = 
    e\ut e$, $\Delta(y_j) = e\ut y_j + y_j\ut e$ and the counit $\eps (e)=1$, 
   $\eps(y_j) = 0$. Let $A= k[x_i\; |\; i\in \Z]/(x_ix_j\; |\; i,j\in \Z)$. Define a linear map $\psi:C\ot A\to A\ot C$ by
$$
\psi (c\ot 1) =1\ot c, \quad \psi(e\ut a) = a\ut e, \quad \psi(y_j\ut x_i) = x_{i+1}
\ut y_{j+1},
$$
  for all $a\in A$, $c\in C$. 
Then $(A,C)_{\psi}$ is an entwining structure, and hence, for all $r,s,p,t\in k$,  $W: A\ot A\to A\ot A$, $X:A\ot C\to A\ot C$, $Z:C\ot C\to C\ot C$, defined by, 
$$
W(1\ot 1) = r1\ot 1, \quad W(x_i\ot 1) = sx_i\ot 1 + (r-s)1\ot x_i, 
$$
$$
W(1 \ot x_i) =  r 1 \ot x_i, \quad W(x_i\ot x_j) = -s x_j\ot x_i,
$$
$$
X(1\ot e) = 1\ot e, \quad X(1\ot y_i) = 1\ot y_i, \quad X(x_i\ut e) = x_i\ut e, \quad X(x_i \ot y_j) = x_{i+1}
\ut y_{j+1}$$
$$
Z(e\ot e) = te\ot e, \quad Z(e\ot y_j) = te\ot y_j + (t-p)y_j\ot e,
$$
$$
 Z(y_j\ot e) = py_j\ot e,\quad  Z(y_i\ot y_j) =- p y_j\ot y_i,
$$
form a Yang-Baxter system.
\end{example}

The final example we consider will illustrate how, in view of Theorem~\ref{thm.main} the knowledge of Yang-Baxter systems can lead to construction (and classification) of entwining structures with invertible entwining map. We use
the full classification of WXZ-systems
 in dimension two, i.e., dim$_k V$= dim$_k V'$=2, obtained in \cite{HlaSno:sol}.

\begin{example}
Guided by the classification of  WXZ-systems in \cite{HlaSno:sol}, for any $s\in k$ ($s\neq -1$), consider  a two-dimensional  algebra $A$ with basis $\{1,x\}$  and with relation  $x^{2} = \frac{1}{s+1}$,  and a two-dimensional coalgebra $C$ with basis
$ \{ e, f \}$ and the coproduct and counit
$$ \Delta (e) = e \ot e + \frac{1}{s+1} f \ot f , \;
 \Delta (f) = e \ot f + f \ot e,\; 
\eps (e) = 1, \;   \eps (f)=0. $$
 Take $W$ and $Z$ as in Theorem~\ref{thm.main} with $r=t=1$ and $p=s$.  Thus explicitly, the operator $W$ reads
$$ W( 1 \ot 1) = 1 \ot 1, \qquad W( 1 \ot x) = 1 \ot x,$$
$$
W( x \ot 1) = (1-s) 1 \ot x + s x \ot 1, \qquad W( x \ot x) = 1 \ot 1 - s x \ot x, $$
or in the matrix form
\begin{equation} \label{w}
W= \begin{pmatrix}
1 & 0 & 0 & 0\\
0 & 1 & 0 & 0\\
0 & 1-s & s & 0\\
1 & 0 & 0 & -s
\end{pmatrix},
\end{equation}
and 
$Z$ is described by a matrix that is the transpose of matrix $W$.

With these operators at hand we look at  the list of solutions for Yang-Baxter systems in
\cite{HlaSno:sol}, to find out that the only relevant cases, i.e., those in which $W$ and $Z$ come from algebras and coalgebras, are the cases 1, 2, 56 and 59 in there.
In the case 1, the entwining map comes out as an ordinary twist. Furthermore, 2 is a special case of  56, the latter described by the matrix
\begin{equation} \label{x}
X= \begin{pmatrix}
1 & 0 & 0 & 0\\
0 & 1 & 0 & 0\\
0 & q & 1 & 0\\
0 & 0 & 0 & -1
\end{pmatrix}.
\end{equation}
The corresponding entwining map $ \psi : C \ot A \rightarrow A \ot C $ comes out as 
\begin{equation}\label{onlypsi}
 \psi (e \ot x) = q 1 \ot f + x \ot e ,\qquad  \psi (f \ot x) = - x \ot f.
\end{equation}

Finally  the case 59 is described by the matrix
\begin{equation} \label{x1}
X=\begin{pmatrix}
s^{-1} & 0 & 0 & q(1-s)\\
0 & 1 & 0 & 0\\
0 & q & s^{-1} & 0\\
0 & 0 & 0 & -1
\end{pmatrix}.
\end{equation}
This means in particular that $ X (1 \ot e) = s^{-1} 1 \ot e + q (1-s) x \ot f $. Thus to satisfy assumptions of Theorem~\ref{thm.main} one has to set $s=1$, in which case the corresponding entwining map coincides with the map given by equations (\ref{onlypsi}). 
\end{example}

We have already explained a connection between Yang-Baxter systems and Yang-Baxter operators obtained by the gluing procedure described in \cite{MajMar:glu}. Given an entwining structure $(A,C)_\psi$ such that $\psi$ is bijective and provided that an appropriate choice of the parameters in Theorem~\ref{thm.main} is made, in addition to the Yang-Baxter operator on $A\oplus C$ discussed earlier, one can construct another Yang-Baxter operator $R$ such that $(R+q^{-1}\id)\circ (R-q\id) = 0$ for a non-zero $q\in k$. Such a Yang-Baxter operator is known as a {\em $q$-Hecke operator}. This follows from  \cite[Example 2.11]{MajMar:glu} and can be stated as the following
\begin{proposition}
Let $(A,C)_\psi$ be an entwining structure and let $W$ and $Z$ be as in Theorem~\ref{thm.main}. Suppose that  $r=t=q$ and
$ s=p=q^{-1} $, for a non-zero $q\in k$, and  that $ \psi $ is bijective. 
Let $ R = W \circ \tau_{A,A} $ and $ R' = Z \circ \tau_{C,C} $ and  define a linear map
$$
 R \oplus_{\psi} R': (A \oplus C) \ot (A \oplus C) \to  (A \oplus C) \ot (A \oplus C)
$$ 
by
$ R \oplus_{\psi} R'|_{A \ot A} = R $,
$ R \oplus_{\psi} R'|_{C\ot C} = R' $, and, for all $ a \in A $, $ c \in C $,
$$ 
(R \oplus_{\psi} R')(a \ot c ) = \psi^{-1}(a \ot c) \qquad 
(R \oplus_{\psi} R')(c \ot a) = \psi(c \ot a)
+ (q-q^{-1})c\ot a.
$$
Then $ R \oplus_{\psi} R'$ is a $q$-Hecke operator.
\end{proposition}
\begin{proof}
Note that $R = R^A_{q,q^{-1}}$, where $R^A_{r,s}$ is  given in equation (\ref{ra}), while $R'=R_C^{q^{-1},q}$, where $R_C^{p,t}$ is given in (\ref{rc}). The minimal polynomials of $R^A_{r,s}$ and  $R_C^{p,t}$ are  $(R^A_{r,s}+s)(R^A_{r,s}-r) $ and $(R_C^{p,t}+p)(R_C^{p,t}-t) $, so that both $R:A\ot A\to A\ot A$ and $R':C\ot C\to C\ot C$ are $q$-Hecke operators. In view of Theorem~\ref{thm.main},  \cite[Example 2.11]{MajMar:glu} implies that also $ R \oplus_{\psi} R'$ is a $q$-Hecke operator.
\end{proof}

\section*{Acknowledgments} 
Tomasz Brzezi\'nski  would like to thank the Engineering and Physical Sciences Research 
Council for an Advanced Fellowship. Florin F.\ Nichita would like to
thank for a Marie Curie Fellowship HPMF-CT-2002-01782.

\end{document}